\documentclass{article}
\usepackage{pamphlet}
\title{A Unified Approach to Algebraic Set Theory}
\author{Benno van den Berg \footnote{Technische Universit\"at Darmstadt, Fachbereich Mathematik, Schlossgartenstr. 7, 64289 Darmstadt, Germany.
berg@mathematik.tu-darmstadt.de.} \and Ieke Moerdijk\footnote{Mathematical Institute, Utrecht University, PO Box 80.010, 3508 TA Utrecht, The Netherlands. moerdijk@math.uu.nl.}}
\date{June 22, 2007}

\begin{document}

\maketitle

\section{Introduction}

This short paper provides a summary of the tutorial on categorical logic given by the second named author at the Logic Colloquium in Nijmegen. Before we go into the subject matter, we would like to express our thanks to the organisers for an excellent conference, and for offering us the opportunity to present this material.

Categorical logic studies the relation between category theory and logical languages, and provides a very efficient framework in which to treat the syntax and the model theory on an equal footing. For a given theory $T$ formulated in a suitable language, both the theory itself and its models can be viewed as categories with structure, and the fact that the models are models of the theory corresponds to the existence of canonical functors between these categories. This applies to ordinary models of first order theories, but also to more complicated topological models, forcing models, realisability and dialectica interpretations of intuitionistic arithmetic, domain-theoretic models of the $\lambda$-calculus, and so on. One of the best worked out examples is that where $T$ extends the theory {\bf HHA} of higher order Heyting arithmetic \cite{lambekscott86}, which is closely related to the Lawvere-Tierney theory of elementary toposes. Indeed, every elementary topos (always taken with a natural numbers object here) provides a categorical model for {\bf HHA}, and the theory {\bf HHA} itself also corresponds to a particular topos, the ``free'' one, in which the true sentences are the provable ones.

The logic of many particular toposes shares features of independence results in set theory. For example, there are very natural constructions of toposes which model {\bf HHA} plus classical logic in which the axiom of choice fails, or in which the continuum hypothesis is refuted. In addition, one easily finds topological sheaf toposes which model famous consistency results of intuitionistic logic, such as the consistency of {\bf HHA} plus the continuity of all real-valued functions on the unit interval, and realisability toposes validating {\bf HHA} plus ``Church's thesis'' (all functions from the natural numbers to itself are recursive). It took some effort (by Freyd, Fourman, McCarthy, Blass and Scedrov \cite{freyd80,fourman80,blassscedrov89,blassscedrov92} and many others), however, to modify the constructions so as to provide models proving the consistency of such statements with {\bf HHA} replaced by an appropriate set theory such as {\bf ZF} or its intuitionistic counterpart {\bf IZF}. This modification heavily depended on the fact that the toposes in question, namely various so-called Grothendieck toposes and Hyland's effective topos \cite{hyland82}, were in some sense defined in terms of sets.

The original purpose of ``algebraic set theory'' \cite{joyalmoerdijk95} was to identify a categorical structure independently of sets, which would allow one to construct models of set theories like {\bf (I)ZF}. These categorical structures were pairs $(\ct{E},\smallmap{S})$ where \ct{E} is a category much like a topos, and \smallmap{S} is a class of arrows in \ct{E} satisfying suitable axioms, and referred to as the class of ``small maps''. It was shown in loc.~cit.~that any such structure gave rise to a model of {\bf (I)ZF}. An important feature of the axiomisation in terms of such pairs $(\ct{E},\smallmap{S})$ is that it is  preserved under the construction of categories of sheaves and of realisability categories, so that the model constructions referred to above become special cases of a general and ``elementary'' preservation result.

In recent years, there has been a lot of activity in the field of algebraic set theory, which is well documented on the web site www.phil.cmu.edu/projects/ast. Several variations and extensions of the the original Joyal-Moerdijk axiomatisation have been developed. In particular, Alex Simpson \cite{simpson99} developed an axiomatisation in which \ct{E} is far from a topos (in his set-up, \ct{E} is not exact, and is only assumed to be a regular category). This allowed him to include the example of classes in {\bf IZF}, and to prove completeness for {\bf IZF} of models constructed from his categorical pairs $(\ct{E},\smallmap{S})$. This approach has been further developed by Awodey, Butz, Simpson and Streicher in their paper \cite{awodeyetal04}, in which they prove a categorical completeness theorem characterising the category of small objects in such a pair $(\ct{E},\smallmap{S})$ (cf.~\reftheo{catcomplforothcat} below), and identify a weak ``basic'' intuitionistic set theory {\bf BIST} corresponding to the core of the categorical axioms in their setting.

In other papers, a variant has been developed which is adequate for constructing models of \emph{predicative} set theories like Aczel's theory {\bf CZF} \cite{aczel78,aczelrathjen01}. The most important feature of this variant is that in the structure $(\ct{E}, \smallmap{S})$, the existence of suitable power objects is replaced by that of inductive W-types. These W-types enabled Moerdijk and Palmgren in \cite{moerdijkpalmgren02} to prove the existence of a model $V$ for {\bf CZF} out of such a structure $(\ct{E}, \smallmap{S})$ on the basis of some exactness assumptions on \ct{E}, and to derive the preservation of (a slight extension of) the axioms under the construction of sheaf categories. This result was later improved by Van den Berg \cite{berg05a}. It is precisely at this point, however, that we believe our current set-up to be superior to the ones in \cite{moerdijkpalmgren02} and \cite{berg05a}, and we will come back to this in some detail in Section 6 below. We should mention here that sheaf models for {\bf CZF} have also been considered by Gambino \cite{gambino06} and to some extent go back to Grayson \cite{grayson83}. Categorical pairs $(\ct{E},\smallmap{S})$ for weak predicative set theories have also been considered by Awodey-Warren \cite{awodeywarren05} and Simpson \cite{simpson05}. (Note, however, that these authors do not consider W-types and only deal with set theories weaker than Aczel's {\bf CZF}.)

The purpose of this paper is to outline an axiomatisation of algebraic set theory which combines the good features of all the approaches mentioned above. More precisely, we will present axioms for pairs $(\ct{E},\smallmap{S})$ which
\begin{itemize}
\item imply the existence in $\ct{E}$ of a universe $V$, which models a suitable set theory (such as {\bf IZF}) (cf.~\reftheo{existencepsfixpoints} below);
\item allow one to prove completeness theorems of the kind in \cite{simpson99} and \cite{awodeyetal04} (cf.~\reftheo{catcompl} and \reftheo{catcomplforothcat} below);
\item work equally well in the predicative context (to construct models of {\bf CZF});
\item are preserved under the construction of sheaf categories, so that the usual topological techniques automatically yield consistency results for {\bf IZF}, {\bf CZF} and similar theories;  
\item hold for realisability categories (cf.~Examples 5.3 and 5.4 and Theorem \reftheo{stabunderrealisability}).
\end{itemize}

Before we do so, however, we will recall the axioms of the systems {\bf IZF} and {\bf CZF} of set theory. In the next Section, we will then present our axioms for small maps, and compare them (in Subsection 3.4) to those in the literature. One of the main features of our axiomatisation is that we do not require the category $\ct{E}$ to be exact, but only to possess quotients of ``small'' equivalence relations. This restricted exactness axiom is consistent with the fact that every object is separated (in the sense of having a small diagonal), and is much easier to deal with in many contexts, in particular those of sheaves. Moreover, together with the Collection axiom this weakened form of exactness suffices for many crucial constructions, such as that of the model $V$ of set theory from the universal small map $E \rightarrow U$, or of the associated sheaf of a given presheaf. In Section 4 we will describe the models of set theory obtained from pairs $(\ct{E},\smallmap{S})$ satisfying our axioms, while Section 5 discusses some examples. Finally in Sections 6 and 7, we will discuss in some detail the preservation of the axioms under the construction of sheaf and realisability categories.

Like the tutorial given at the conference, this exposition is necessarily concise, and most of the proofs have been omitted. With the exception perhaps of Sections 6 and 7, these proofs are often suitable adaptations of existing proofs in the literature, notably \cite{joyalmoerdijk95,simpson99,moerdijkpalmgren02,awodeyetal04,berg05}. A complete exposition with full proofs will appear as \cite{bergmoerdijk07}.

We would like to thank Thomas Streicher, Jaap van Oosten and the anonymous referees for their comments on an earlier draft of this paper, and Thomas Streicher in particular for suggesting the notion of a display map defined in Section 7.

\section{Constructive set theories}

In this Section we recall the axioms for the two most prominent constructive variants of Zermelo-Fraenkel set theory, {\bf IZF} and {\bf CZF}. Like ordinary {\bf ZF}, these two theories are formulated in first-order logic with one non-logical symbol $\epsilon$. But unlike ordinary set theory, these theories are constructive, in that their underlying logic is intuitionistic.

In the formulation of the axioms, we use the following standard abbreviations:
$ \exists x \epsilon a \, (\ldots)$ for $ \exists x \, (x \epsilon a \land \ldots)$, and $\forall x \epsilon a \, (\ldots)$ for $ \forall x \, (x \epsilon a \rightarrow \ldots)$. Recall also that a formula is called \emph{bounded}, when all the quantifiers it contains are of one of these two forms. Finally, a formula of the form $\forall x \epsilon a \, \exists y \epsilon b \, \phi \land \forall y \epsilon b \, \exists x \epsilon a \, \phi$ will be abbreviated as:
\[ \mbox{B}(x \epsilon a, y \epsilon b) \, \phi. \]

The axioms which both theories have in common are (the universal closures of):
\begin{description}
\item[Extensionality:] $\forall x \, ( \, x \epsilon a \leftrightarrow x \epsilon b \, ) \rightarrow a = b$.
\item[Empty set:] $\exists x  \, \forall y  \, \lnot y \epsilon x $.
\item[Pairing:] $\exists x \, \forall y \, (\,  y \epsilon x \leftrightarrow y = a \lor y = b \, )$.
\item[Union:] $\exists x  \, \forall y \, ( \, y \epsilon x \leftrightarrow \exists z \epsilon a \, y \epsilon z  \, )$.
\item[$\epsilon$-induction:] $\forall x \, (\forall y  \epsilon x \, \phi(y) \rightarrow \phi(x)) \rightarrow \forall x \, \phi(x)$
\item[Bounded separation:] $\exists x \,  \forall y \, ( \, y \epsilon x \leftrightarrow y \epsilon a \land \phi(y) \, ) $, for any bounded formula $\phi$ in which $a$ does not occur.
\item[Strong collection:] $\forall x \epsilon a \, \exists y \, \phi(x,y) \rightarrow \exists b \, \mbox{B}(x \epsilon a, y \epsilon b) \, \phi$.
\item[Infinity:] $\exists a \, ( \, \exists  x \, x \epsilon a \, ) \land ( \, \forall x  \epsilon a \, \exists y \epsilon a \, x \epsilon y \, )$.
\end{description}

One can obtain an axiomatisation for the constructive set theory {\bf IZF} by adding to the axioms above the following two statements:
\begin{description}
\item[Full separation:] $\exists x \,  \forall y \, ( \, y \epsilon x \leftrightarrow y \epsilon a \land \phi(y) \, ) $, for any formula $\phi$ in which $a$ does not occur.
\item[Power set axiom:] $\exists x \, \forall y \, ( \, y \epsilon x \leftrightarrow y \subseteq a \, )$.
\end{description}

To obtain the predicative constructive set theory {\bf CZF},  one should add instead the following axiom (which is a weakening of the Power Set Axiom):
\begin{description}
\item[Subset collection:] $\exists c \, \forall z \, ( \forall x \epsilon a \, \exists y \epsilon b \, \phi(x,y,z) \rightarrow \exists d \epsilon c \, \mbox{B}(x \epsilon a, y \epsilon d) \, \phi(x, y, z)) $.
\end{description}
The Subset Collection Axiom has a more palatable formulation (equivalent to it relative to the other axioms), called Fullness (see \cite{aczelrathjen01}). Write ${\bf mv}(a, b)$ for the class of all multi-valued functions from a set $a$ to a set $b$, i.e.~relations $R$ such that $\forall x \epsilon a \, \exists y \epsilon b \, (x,y) \epsilon R$.
\begin{description}
\item[Fullness:] $\exists u \, (u \subseteq {\bf mv}(a,b) \land \forall v \epsilon {\bf mv}(a,b) \, \exists w \epsilon u \, (w \subseteq v))$.
\end{description}
Using this formulation, it is also easier to see that Subset Collection implies Exponentiation, the statement that the functions from a set $a$ to a set $b$ form a set.

\section{Categories with small maps}

Here we introduce the categorical structure which is necessary to model set theory. The structure is that of a category \ct{E} equipped with a class of morphisms \smallmap{S}, satisfying certain axioms and being referred to as the \emph{class of small maps}. The canonical example is the one where \ct{E} is the category of classes in a model of some weak set theory, and morphisms between classes are small in case all the \emph{fibres} are sets. More examples will follow in Section 5. In Section 4, we will show that these axioms actually provide us with the means of constructing models of set theory.

\subsection{Axioms}

In our work, the underlying category \ct{E} is a Heyting category with sums. More precisely, \ct{E} satisfies the following axioms (for an excellent account of the notions involved, see \cite[Part A1]{johnstone02a}):
\begin{itemize}
\item \ct{E} is cartesian, i.e.~it has finite limits.
\item \ct{E} is regular, i.e.~every morphism factors as a cover followed by a mono and covers are stable under pullback.
\item \ct{E} has finite disjoint and stable coproducts.
\item \ct{E} is Heyting, i.e.~for any morphism \func{f}{X}{Y} the functor \func{f^*}{{\rm Sub} (Y)}{{\rm Sub} (X)} has a right adjoint $\forall_f$.
\end{itemize}
This expresses precisely that \ct{E} is a categorical structure suitable for modelling a typed version of first-order intuitionistic logic with finite product and sum types.

We now list the axioms that we require to hold for a class of small maps, extending the axioms for a class of open maps (see \cite{joyalmoerdijk95}). We will comment on the relation between our axiomatisation and existing alternatives in Section 3.4 below.

The axioms for a class of open maps \smallmap{S} are:
\begin{description}
\item[(A1)] (Pullback stability) In any pullback square
\diag{ D \ar[d]_g \ar[r] & B \ar[d]^f \\
C \ar[r]_p & A, }
where $f \in \smallmap{S}$, also $g \in \smallmap{S}$.
\item[(A2)] (Descent) Whenever in a pullback square as above, $g \in \smallmap{S}$ and $p$ is a cover, $f \in \smallmap{S}$.
\item[(A3)] (Sums) Whenever $X \rTo Y$ and $X'\rTo Y'$ belong to \smallmap{S}, so does $X + X' \rTo Y + Y'$.
\item[(A4)] (Finiteness) The maps $0 \rTo 1$, $1 \rTo 1$ and $2 = 1+1 \rTo 1$ belong to \smallmap{S}.
\item[(A5)] (Composition) \smallmap{S} is closed under composition.
\item[(A6)] (Quotients) In any commutative triangle
\diag{ Z \ar@{->>}[rr]^p \ar[dr]_g & & Y \ar[dl]^f \\
& X, & }
where $p$ is a cover and $g$ belongs to \smallmap{S}, so does $f$.
\end{description}
These axioms are of two kinds: the axioms {\bf (A1-3)} express that the property we are interested in is one of the \emph{fibres} of maps in \smallmap{S}. The others are more set-theoretic: {\bf (A4)} says that the collections containing 0, 1 or 2 elements are sets. {\bf (A5)} is a union axiom: the union of a small disjoint family of sets is again a set. Finally, {\bf (A6)} is a form of replacement: the image of a set is again a set.

We will always assume that a class of small maps \smallmap{S} satisfies the following two additional axioms, familiar from \cite{joyalmoerdijk95}:
\begin{description}
\item[(C)] (Collection) Any two arrows \func{p}{Y}{X} and \func{f}{X}{A} where $p$ is a cover and $f$ belongs to \smallmap{S} fits into a quasi-pullback diagram\footnote{Recall that a commutative square in a regular category is called a quasi-pullback if the unique arrow from the initial vertex of the square to the inscribed pullback is a cover.} of the form
\diag{ Z \ar[d]_g \ar[r] & Y \ar@{->>}[r]^p & X \ar[d]^f \\
B \ar@{->>}[rr]_h & & A,}
where $h$ is a cover and $g$ belongs to \smallmap{S}.
\item[(R)] (Representability, see \refrema{strongrepres}) There exists a small map \func{\pi}{E}{U} (a ``universal small map'') such that for every small map \func{f}{X}{Y} there is a diagram of the shape
\diag{X \ar[d]_f & A \ar[d] \ar[r] \ar@{->>}[l] & E \ar[d]^{\pi} \\
Y & B \ar[r] \ar@{->>}[l]^p & U, }
where the left square is a quasi-pullback, the right square is a pullback and $p$ is a cover.
\end{description}
The collection principle {\bf (C)} expresses that in the internal logic it holds that for any cover \func{p}{Y}{X} with small codomain there is a cover $Z \rTo X$ with small domain that factors through $p$, while {\bf (R)} says that there is a (necessarily class-sized) family of sets $(E_u)_{u \in U}$ such that any set is covered by one in this family.

The next requirement is also part of the axioms in \cite{joyalmoerdijk95}. For a morphism \func{f}{X}{Y}, the pullback functor \func{f^*}{\ct{E}/Y}{\ct{E}/X} always has a left adjoint $\Sigma_f$ given by composition\footnote{We will write $X^*$ and $\Sigma_X$ for $f^*$ and $\Sigma_f$, where $f$ is the unique map $X \rTo 1$.}. It has a right adjoint $\Pi_f$ only when $f$ is exponentiable.
\begin{description}
\item[($\Pi$E)] (Existence of $\Pi$) The right adjoint $\Pi_f$ exists, whenever $f$ belongs to \smallmap{S}.
\end{description}
This intuitively means that for any set $A$ and class $X$ there is a class of functions from $A$ to $X$.

When $f$ is exponentiable, one can define an endofunctor $P_f$ (the polynomial functor associated with $f$) as the composition:
\[ P_f = \Sigma_Y \Pi_f X^*. \]
Its initial algebra (whenever it exists) is called the W-type associated to $f$. For extensive discussion and examples of these W-types we refer the reader to \cite{moerdijkpalmgren00,berg05,gambinohyland04}. We impose the axiom (familiar from \cite{moerdijkpalmgren00,gambino05}):
\begin{description}
\item[(WE)] (Existence of $W$) The W-type associated to any map \func{f}{X}{Y} in \smallmap{S} exists.
\end{description}
In non-categorical terms this means that for a signature consisting of a (possibly class-sized) number of term constructors each of which has an arity forming a set, the free term algebra exists (but maybe not as a set).

The following two axioms are necessary to have bounded separation as an internally valid principle (see \refrema{boundedsepinbasicset}). For this purpose we need a piece of terminology: call a subobject \diag{m: A \ar@{ >->}[r] & X} \smallmap{S}-bounded, whenever $m$ belongs to \smallmap{S}; note that the \smallmap{S}-bounded subobjects form a submeetsemilattice of Sub($X$). We impose the following axiom:
\begin{description}
\item[(HB)] (Heyting axiom for bounded subobjects) For any small map \func{f}{Y}{X} the functor \func{\forall_f}{{\rm Sub}(Y)}{{\rm Sub}(X)} maps \smallmap{S}-bounded subobjects to \smallmap{S}-bounded subobjects.
\end{description}
In addition, we require that all equalities are bounded. Call an object $X$ separated, when the diagonal $\Delta: X \rTo X \times X$ is small. We furthermore impose (see \cite{awodeyetal04}):
\begin{description}
\item[(US)] (Universal separation) All objects are separated.
\end{description}

We finally demand a limited form of \emph{exactness}, by requiring the existence of quotients for a restricted class of equivalence relations. To formulate this categorically, we recall the following definitions. Two parallel arrows
\diag{
  R \ar@<.8ex>[r]^{r_0} \ar@<-.8ex>[r]_-{r_1} & X
}
in category \ct{E} form an \emph{equivalence relation} when for any object $A$ in \ct{E} the induced function
\[ \mbox{Hom}(A,R) \rTo \mbox{Hom}(A, X) \times \mbox{Hom}(A, X) \]
is an injection defining an equivalence relation on the set $\mbox{Hom}(A,X)$. We call an equivalence relation bounded, when $R$ is a bounded subobject of $X \times X$. A morphism \func{q}{X}{Q} is called the \emph{quotient} of the equivalence relation, if the diagram
\diag{
  R \ar@<.8ex>[r]^{r_0} \ar@<-.8ex>[r]_-{r_1} & X \ar[r]^q & Q
}
is both a pullback and a coequaliser. In this case, the diagram is called \emph{exact}. The diagram is called \emph{stably exact}, when for any \func{p}{P}{Q} the diagram
\diag{
  p^{*}R \ar@<.8ex>[r]^{p^{*}r_0} \ar@<-.8ex>[r]_-{p^{*}r_1} & p^{*}X \ar[r]^{p^{*}q} & p^{*}Q
}
is also exact. If the quotient completes the equivalence relation to a stably exact diagram, we call the quotient stable.

In the presence of {\bf (US)}, any equivalence relation that has a (stable) quotient, must be bounded. So our last axiom imposes the maximum amount of exactness that can be demanded:
\begin{description}
\item[(BE)] (Bounded exactness) All $\smallmap{S}$-bounded equivalence relations have stable quotients.
\end{description}

This completes our definition of a class of small maps. A pair $(\ct{E}, \smallmap{S})$ satisfying the above axioms now will be called \emph{a category with small maps}. 

When a class of small maps \smallmap{S} has been fixed, we call a map $f$ small if it belongs to \smallmap{S}, an object $A$ small if $A \rTo 1$ is small, a subobject \func{m}{A}{X} small if $A$ is small, and a relation $R \subseteq C \times D$ small if the composite
\[ R \subseteq C \times D \rTo D \]
is small.

We conclude this Subsection with some remarks on a form of exact completion relative to a class of small maps. As a motivation, notice that axiom {\bf (BE)} is not satisfied in our canonical example, where \ct{E} is the category of classes in a model of some weak set theory. To circumvent this problem, we will prove the following theorem in our companion paper \cite{bergmoerdijk07}:
\begin{theo}{consofBE}
The axiom {\bf (BE)} is conservative over the other axioms, in the following precise sense. Any category \ct{E} equipped with a class of maps \smallmap{S} satisfying all axioms for a class of small maps except {\bf (BE)} can be embedded in a category $\overline{\ct{E}}$ equipped with a class of small maps $\overline{\smallmap{S}}$ satisfying all the axioms, including {\bf (BE)}. Moreover, the embedding \func{\bf y}{\ct{E}}{\overline{\ct{E}}} is fully faithful, bijective on subobjects and preserves the structure of a Heyting category with sums, hence preserves and reflects validity of statements in the internal logic. Finally, it also preserves and reflects smallness, in the sense that ${\bf y}f$ belongs to $\overline{\smallmap{S}}$ iff $f$ belongs to $\smallmap{S}$.
\end{theo}
The category $\overline{\ct{E}}$ is obtained by formally adjoining quotients for bounded equivalence relations, as in \cite{carboniceliamagno82,carboni95}. Furthermore, a map \func{g}{B}{A} in $\overline{\ct{E}}$ belongs to $\overline{\smallmap{S}}$ iff it fits into a quasi-pullback square
\diag{{\bf y}D \ar[d]_{{\bf y}f} \ar@{->>}[r] & B \ar[d]^g \\
{\bf y}C \ar@{->>}[r] & A, }
with $f$ belonging to \smallmap{S} in \ct{E}.

\subsection{Consequences}

Among the consequences of these axioms we list the following.

\begin{rema}{stabilityunderslicing} 
For any object $X$ in \ct{E}, the slice category $\ct{E}/X$ is equipped with a class of small maps $\smallmap{S}/X$, by declaring that an arrow $p \in \ct{E}/X$ belongs to $\smallmap{S}/X$ whenever $\Sigma_X f$ belongs to $\smallmap{S}$. Any further requirement for a class of small maps should be stable under slicing in this sense, if it is to be a sensible addition. We will not explicitly check this every time we introduce a new axiom, and leave this to the reader.
\end{rema}

\begin{rema}{boundedsepinbasicset} In a category \ct{E} with small maps the following internal form of ``bounded separation'' holds. If $\phi(x)$ is a formula in the internal logic of \ct{E} with free variable $x \in X$, all whose basic predicates are bounded, and contains existential and universal quantifications $\exists_f$ and $\forall_f$ only along small maps $f$, then
\[ A = \{ x \in X \, | \, \phi(x) \} \subseteq X \]
defines a bounded subobject of $X$. In particular, smallness of $X$ implies smallness of $A$.
\end{rema}

\begin{rema}{strongrepres}
It follows from the axioms that any class of small maps \smallmap{S} is also representable in the stronger sense that there is a universal small map \func{\pi}{E}{U} such that for every small map \func{f}{X}{Y} there is a diagram of the shape
\diag{X \ar[d]_f & A \ar[d] \ar[r] \ar@{->>}[l] & E \ar[d]^{\pi} \\
Y & B \ar[r] \ar@{->>}[l]^p & U, }
where the left square is a pullback, the right square is a pullback and $p$ is a cover. Actually, this is how representability was stated in \cite{joyalmoerdijk95}. We have chosen the weaker formulation {\bf (R)}, because it is easier to check in some examples.
\end{rema}

\begin{rema}{smallmapspowerful} Using the axioms {\bf ($\Pi$E), (R), (HB)} and {\bf (BE)}, it can be shown along the lines of Theorem 3.1 in \cite{joyalmoerdijk95} that for any class of small maps the following axiom holds:
\begin{description}
\item[(PE)] (Existence of power class functor) For any object $C$ in \ct{E} there exists a power object $\spower C$ and a small relation $\in_C \subseteq C \times \spower C$ such that, for any $D$ and any small relation $R \subseteq C \times D$, there exists a unique map \func{\rho}{D}{\spower C} such that the square:
\diag{ R \ar@{ >->}[d] \ar[r] & \in_C \ar@{ >->}[d] \\
C \times D \ar[r]_{1 \times \rho} & C \times \spower C}
is a pullback.
\end{description}
In addition, one can show that the object $\spower C$ is unique (up to isomorphism) with this property, and that the assignment $C \mapsto \spower C$ is functorial. 

A special role is played by $\Omega_b = \spower 1$, what one might call the object of bounded truth-values, or the bounded subobject classifier. There are a couple of observations one can make: bounded truth-values are closed under small infima and suprema, implication, and truth and falsity are bounded truth-values. A subobject \func{m}{A}{X} is bounded, when the assertion ``$x \in A$'' has a bounded truth-value for any $x \in X$, as such bounded subobjects are classified by maps $X \rTo \Omega_b$.
\end{rema}

\begin{rema}{smallthingspretopos} (See \cite{awodeywarren05}.) When \ct{E} is a category with a class of small maps \smallmap{S}, and we fix an object $X \in \ct{E}$, we can define a full subcategory \smallmap{S}$_X$ of \ct{E}$/X$, whose objects are small maps into $X$. The category \smallmap{S}$_X$ is a Heyting pretopos, and the inclusion into \ct{E}$/X$ preserves this structure; this was proved in \cite{awodeywarren05}. This result can be regarded as a kind of categorical ``soundness'' theorem, in view of the following corresponding ``completeness'' theorem, which is analogous to Grothendieck's result that every pretopos arises as the coherent objects in a coherent topos (see \cite[Section D.3.3]{johnstone02b}).
\begin{theo}{catcompl}
Any Heyting pretopos \ct{H} arises as the category of small objects $\smallmap{S}_1$ in a category \ct{E} with a class of small maps \smallmap{S}.
\end{theo}
This theorem was proved in \cite{awodeywarren05}, where, following \cite{awodeyetal04}, the objects in \ct{E} were called the \emph{ideals} over \ct{H}.
\end{rema}

\subsection{Strengthenings}

For the purpose of constructing models of important (constructive) set theories, we will consider the following additional properties which a class of small maps may enjoy.

\begin{description}
\item[(NE)] (Existence of nno) The category \ct{E} possesses a natural numbers object.
\item[(NS)] (Smallness of nno) In addition, it is small.
\end{description}
There is no need to impose {\bf (NE)}, as it follows from {\bf (WE)}. The axiom {\bf (NS)} is necessary for modelling set theories with Infinity. The property {\bf (PE)} in \refrema{smallmapspowerful} has a similar strengthening, corresponding to the Power set Axiom:
\begin{description}
\item[(PS)] (Smallness of power classes) For each $X$ the $\spower$-functor on $\ct{E}/X$ preserves smallness of objects over $X$.
\end{description}
Both {\bf (NS)} and {\bf (PS)} were formulated in \cite{joyalmoerdijk95} for the purpose of modelling {\bf IZF}.
\begin{rema}{smallthingstopos} (Cf.~\cite{awodeyetal04}.)
Let $X$ be an object in a category with small maps $(\ct{E}, \smallmap{S})$ satisfying ${\bf (PS)}$. The category \smallmap{S}$_X$ is a topos, and the inclusion into \ct{E}$/X$ preserves this structure. In fact, every topos arises in this way:
\begin{theo}{catcomplforothcat}
Any topos \ct{H} arises as the category of small objects $\smallmap{S}_1$ in a category equipped with a class of small maps satisfying {\bf (PS)}.
\end{theo}
Like \reftheo{catcompl}, this is proved in \cite{awodeyetal04} using the ideal construction. 
\end{rema}

We will also need to consider requirements corresponding to the axioms of Full Separation and Fullness. To Full Separation corresponds the following axiom, introduced in \cite{joyalmoerdijk95}:
\begin{description}
\item[(M)] All monos are small.
\end{description}
A categorical axiom corresponding to Fullness was first stated in \cite{bergdemarchi06}. In order to formulate it, we need to introduce some notation. For two morphisms $A \rTo X$ and $B \rTo X$, we will denote by $M_X(A, B)$  the poset of multi-valued functions from $A$ to $B$ over $X$, i.e.~jointly monic spans in $\ct{E}/X$,
\diag{ A & \ar@{->>}[l] P \ar[r] & B}
with $P \rTo X$ small and the map to $A$ a cover. By pullback, any \func{f}{Y}{X} determines an order preserving function
\[ \func{f^*}{M_X(A, B)}{M_Y(f^*A, f^*B)}. \]
\begin{description}
\item[(F)] For any two small maps $A \rTo X$ and $B \rTo X$, there exist a cover $p: X' \rTo X$, a small map $f: C \rTo X'$ and an element $P \in M_C(f^*p^*A, f^*p^*B)$, such that for any $g: D \rTo X'$ and $Q \in M_D(g^*p^*A, g^*p^*B)$, there are morphisms $x: E \rTo D$ and $y: E \rTo C$, with $gx = fy$ and $x$ a cover, such that $x^* Q \geq y^* P$.
\end{description}
Though complicated, it is ``simply'' the Kripke-Joyal translation of the statement that there is for any pair of small objects $A$ and $B$, a small collection $P$ of multi-valued relations between $A$ and $B$, such that any multi-valued relation contains one in $P$.

\subsection{Relation to other settings}

The axioms for a category with small maps $(\ct{E}, \smallmap{S})$ as we have presented them are very close to the original axioms as presented by Joyal and Moerdijk on pages 6-8 of their book \cite{joyalmoerdijk95}. We only require the weak form of exactness of {\bf (BE)} (instead of ordinary exactness), and added the axioms {\bf (WE), (HB)} and {\bf (US)}.

Since the appearance of \cite{joyalmoerdijk95}, various axiomatisations have been proposed, which can roughly be subdivided into three groups. To the first group belong axiom systems extending the original presentation in \cite{joyalmoerdijk95}. Already in \cite{joyalmoerdijk95}, it is shown how to extend these axioms for the purpose of obtaining models for {\bf IZF}, and this is followed up in \cite{kouwenhovenvanoosten05}. In \cite{gambino05} Gambino introduces an extension of the original axiomatisation leading to models of predicative set theories.

A second group of papers starts with Simpson's \cite{simpson99} and comprises \cite{awodeyetal04,butz03,simpson05,awodeyforssell05,awodeywarren05}. In these axiomatisations, the following axioms which are here taken as basic are regarded as optional features: the Collection Axiom {\bf (C)}, Bounded Exactness {\bf (BE)}, and also {\bf (WE)} (although they all hold in the category of ideals). Instead, the existence of a \spower-functor as in {\bf (PE)} is postulated, as is a model of set theory, either in the form of a universe, or a universal object. In the approach taken here, these are properties derived from the existence of a universal small map \func{\pi}{E}{U}. Part of the purpose of this paper is to make clear that the results for axiom systems in \cite{simpson99,awodeyetal04,awodeywarren05} also hold for our axiomatisation. We list the achievements in order to make a comparison possible: in \cite{simpson99}, Simpson obtained a set-theoretic completeness result for an impredicative set theory (compare \reftheo{completeness}). Then in \cite{awodeyetal04}, Awodey, Butz, Simpson and Streicher prove a categorical completeness result of which our \reftheo{catcomplforothcat} is variant. A predicative version of this result which does not involve W-types but is otherwise analogous to \reftheo{catcompl} above, was then proved by Awodey and Warren in \cite{awodeywarren05}.

The fact that our set-up contains the Collection Axiom {\bf (C)} makes it less appropriate for modelling set theories based on the axiom of Replacement. However, in our theory this Collection Axiom plays a crucial r\^ole: for example, in the construction of the initial ZF-algebra from W-types (see \reftheo{existencepsfixpoints} below), or in showing the existence of the associated sheaf functor.

A third group of papers starts with \cite{moerdijkpalmgren02}, and continues with \cite{berg05a,berg06}. These axiomatisations have a flavour different from the others, because here the axioms for a class of small maps do not extend the axioms for a class of open maps, as the Quotient Axiom {\bf (A6)} is dropped. The aim of Moerdijk and Palmgren in \cite{moerdijkpalmgren02} was to find an axiomatisation related to Martin-L\"of's predicative type theory which included the category-theoretic notion of a W-type, from which models of Aczel's {\bf CZF} could be constructed. We will point out below that the same is true here (in fact, we can construct models of {\bf CZF} proper, rather than of something less or more). Another concern of \cite{moerdijkpalmgren02}, which is also the topic of Van den Berg's paper \cite{berg05a}, is the stability of the notion of category with small maps under sheaves. The earlier results in \cite{moerdijkpalmgren02} and \cite{berg05a} concerning sheafification were less than fully satisfactory. For the notion of category with small maps explained here, the theory of sheaves can be developed very smoothly (see Section 6), using the combination of the axioms {\bf (BE)} and {\bf (US)}. We consider this one of the main advantages of the present axiomatisation.

\section{Models of set theory}

For the purpose of discussing models of set theory, we recall from \cite{joyalmoerdijk95} the notion of a \emph{ZF-algebra} in a category with small maps $(\ct{E}, \smallmap{S})$. A ZF-algebra $V$ is an object in \ct{E} equipped with two independent algebraic structures: on the one hand, it is an (internal) poset with small (in the sense of \smallmap{S}) sups. On the other hand, it is equipped with an endomap $s: V \rTo V$, called ``successor''. A morphism of ZF-algebras should preserve both these structures: the small suprema, and the successor.

A crucial result is the following:
\begin{theo}{existencepsfixpoints}
In any category with small maps $(\ct{E}, \smallmap{S})$, the initial ZF-algebra exists.
\end{theo}
This theorem can be proved along the lines of \cite{moerdijkpalmgren02}. Indeed, one can consider the W-type associated to the universal small map \func{\pi}{E}{U}. One can then show that the equivalence relation given by bisimulation is bounded so that the quotient exists. This quotient is the initial ZF-algebra (more details will appear in \cite{bergmoerdijk07}). 
This initial ZF-algebra has a natural interpretation as a model of set theory. We think of the order as inclusion, suprema as union, and $sx$ as $\{ x \}$. This suggests to define membership as:
\begin{eqnarray*}
x \, \epsilon \, y & := & sx \leq y.
\end{eqnarray*}
Since \ct{E} is a Heyting category, one can ask oneself the question which set-theoretic statements the structure $(V, \epsilon)$ satisfies in the internal logic of \ct{E}. The answer is given by the following theorem, whose second part was proved in \cite{joyalmoerdijk95} (the first part can be proved in a similar manner):
\begin{theo}{soundness} Let $(\ct{E}, \smallmap{S})$ be a category with small maps in which the natural numbers object is small (so {\bf (NS)} holds). 
\begin{enumerate}
\item If $(\ct{E}, \smallmap{S})$ satisfies the Fullness Axiom {\bf (F)}, then the initial ZF-algebra models {\bf CZF}. 
\item If $(\ct{E}, \smallmap{S})$ satisfies the Power Set Axiom {\bf (PS)} and the Separation axiom {\bf (M)}, then the initial ZF-algebra models {\bf IZF}.
\end{enumerate}
\end{theo}

\begin{rema}{classetth}
To obtain models for classical set theories, one may work in Boolean categories. Initial ZF-algebras in such categories validate classical logic, and therefore model classical set theories.
\end{rema}

As a counterpart to \reftheo{soundness} we can formulate a completeness theorem:
\begin{theo}{completeness} The semantics of \reftheo{soundness} is complete for both {\bf CZF} and {\bf IZF} in the following strong sense.
\begin{enumerate}
\item There is a category with small maps $(\ct{E}, \smallmap{S})$ satisfying {\bf (NS)} and {\bf (F)} such that its initial ZF-algebra $V$ has the property that, for any sentence $\phi$ in the language of set theory:
\[ V \models \phi \Leftrightarrow {\rm \bf CZF} \vdash \phi. \]
\item There is a category with small maps $(\ct{E}, \smallmap{S})$ satisfying {\bf (NS)}, {\bf (M)} and {\bf (PS)}, such that its initial ZF-algebra $V$ has the property that, for any sentence $\phi$ in the language of set theory:
\[ V \models \phi \Leftrightarrow {\rm \bf IZF} \vdash \phi. \]
\end{enumerate}
\end{theo}
To prove this theorem one builds the syntactic category of classes and a ZF-algebra $V$ such that validity in $V$ is the same as derivability in the appropriate set theory. Problems concerning {\bf (BE)} are, of course, solved by appealing to \reftheo{consofBE}. The first person to prove a completeness result in this manner was Alex Simpson in \cite{simpson99} for an impredicative set theory. A predicative variation is contained in \cite{awodeywarren05} and \cite{gambino05}. 

\begin{rema}{complrelativetomodel}
Every (ordinary, classical) set-theoretic model $(M, \epsilon)$ is also subsumed in our account, because every such model there is an initial ZF-algebra $V_M$ in a category with small maps $(\ct{E}_M, \smallmap{S}_M)$ having the property that for any set-theoretic sentence $\phi$:
\[ V_M \models \phi \Leftrightarrow M \models \phi. \]
$\ct{E}_M$ is of course the category of classes in the model $M$, with those functional relations belonging to $\smallmap{S}_M$ that the model believes to have sets as fibres, extended using \reftheo{consofBE} so as to satisfy {\bf (BE)}. One could prove completeness of our categorical semantics for classical set theories along these lines. 
\end{rema}

\section{Examples}

We recall from \cite{joyalmoerdijk95} the basic examples of categories satisfying our axioms.

\begin{kopje}{Example} The canonical example is the following. Let \ct{E} be the category of classes in some model of set theory, and declare a morphism \func{f}{X}{Y} to be small, when all its fibres are sets. If the set theory is strong enough, this will satisfy all our axioms, except for {\bf (BE)}, but an appeal to \reftheo{consofBE} will resolve this issue.
\end{kopje}

\begin{kopje}{Example}
Let \ct{E} be a category of sets (relative to some model of ordinary set theory, say), and let $\kappa$ be an infinite regular cardinal. Declare \func{f}{X}{Y} to be small, when all fibres of $f$ have cardinality less than $\kappa$. This will validate all our basic axioms, as well as {\bf (M)}. When $\kappa \gt \omega$, {\bf (NS)} will also hold, and when $\kappa$ is inaccessible, {\bf (PS)} and {\bf (F)} will hold.
\end{kopje}

\begin{kopje}{Example} The following two examples are related to realisability, and define classes of small maps on the effective topos \Eff \, (see \cite{hyland82}). Recall that there is an adjoint pair of functors $\Gamma \ladj \nabla$, where $\Gamma = \Eff(1, -): \Eff \rTo \Sets$ is the global sections functor. Fix a regular cardinal $\kappa \gt \omega$, and declare \func{f}{X}{Y} to be small, whenever there is a quasi-pullback square
\diag{ Q \ar@{->>}[r] \ar[d]_g & X \ar[d]^f \\
P \ar@{->>}[r]_p & Y }
with $p$ a cover, and $g$ a morphism between projectives such that $\Gamma g$ is $\kappa$-small, in the sense of the previous example. This example was further studied by Kouwenhoven and Van Oosten in \cite{kouwenhovenvanoosten05}, and shown to lead to McCarty's realisability model of set theory for an inaccessible cardinal $\kappa$ (see \cite{mccarty84}).
\end{kopje}

\begin{kopje}{Example} Another class of small maps on $\Eff$ is given as follows. Call a map \func{f}{X}{Y} small, whenever the statement that all its fibres are subcountable is true in the internal logic of $\Eff$ (a set is subcountable, when it is the quotient of a subset of the natural numbers). These maps were studied in \cite{hylandrobinsonrosolini90} and dubbed ``quasi-modest'' in \cite{joyalmoerdijk95}. The first author showed they lead to a model of {\bf CZF} in which all sets are subcountable, and therefore refutes the Power Set Axiom (see \cite{berg06}). He also showed the model is the same as the one contained in \cite{streicher05} and \cite{lubarsky06}. 
\end{kopje}

\begin{kopje}{Example} Once again, fix an infinite regular cardinal $\kappa$, and let \ct{C} be a subcanonical site which is $\kappa$-small, in the sense that every covering family has cardinality strictly less than $\kappa$. We say that a sheaf $X$ is $\kappa$-small, whenever it is covered by a collection of representables whose cardinality is less than $\kappa$. Finally, a morphism \func{f}{X}{Y} will be considered to be $\kappa$-small, whenever for any map \func{y}{C}{Y} from a representable $C \in \ct{C}$ the pullback $f^{-1}(y)$ as in
\diag{ f^{-1}(y) \ar[r] \ar[d] & X \ar[d]^f \\
C \ar[r]_y & Y }
is a $\kappa$-small sheaf. This can again be shown to satisfy all our basic axioms. Also, when $\kappa \gt \omega$, {\bf (NS)} will hold, and so will {\bf (PS)} and {\bf (F)}, when $\kappa$ is inaccessible.
\end{kopje}

\section{Predicative sheaf theory}

The final example of the previous Section, that of sheaves, can be internalised, in a suitable sense. Starting from a category with small maps $(\ct{E}, \smallmap{S})$, and an appropriate site \ct{C} in \ct{E}, one can build the category ${\rm Sh}_{\ct{E}}(\ct{C})$ of internal sheaves over \ct{C}, which is again a Heyting category with stable, disjoint sums. Furthermore, there is a notion of small maps between sheaves, turning it into a category with small maps. In fact, stability of our notion of a category with small maps under sheaves is one of its main assets. Here we will limit ourselves to formulating precise statements, leaving the proofs for \cite{bergmoerdijk07}.

For the site \ct{C} we assume first of all that the underlying category is small, in that the object of objects $C_0$ and of arrows $C_1$ are both small. By a \emph{sieve} on $a \in C_0$ we mean a \emph{small} collection of arrows into $a$ closed under precomposition. We assume that the collection of covering sieves ${\rm Cov}(a)$ on an object $a \in C_0$ satisfies the following axioms:
\begin{description}
\item[(M)] The maximal sieve $M_a = \{ f \in C_1 \, | \, \mbox{cod}(f) = a \}$ belongs to ${\rm Cov}(a)$.
\item[(L)] For any $U \in {\rm Cov}(a)$ and morphism \func{f}{b}{a}, the sieve $f^*U = \{ \func{g}{c}{b} \, | \, fg \in U \}$ belongs to ${\rm Cov}(b)$.
\item[(T)] If $T$ is a sieve on $a$, such that for a fixed $U \in {\rm Cov}(a)$ any pullback $h^*T$ along a map $\func{h}{b}{a} \in U$ is an element of ${\rm Cov}(b)$, then $T \in {\rm Cov}(a)$.
\end{description}
The definition of an (internal) presheaf and sheaf is as usual.

Using the bounded exactness of $(\ct{E}, \smallmap{S})$ and assuming that the relation $S \in {\rm Cov}(a)$ is bounded, one can show the existence of the associated sheaf functor (the cartesian left adjoint for the inclusion of sheaves into presheaves). This functor can then be used to prove in the usual way that the sheaves form a Heyting category with stable and disjoint sums.

As the small maps between sheaves we take those that are ``pointwise small''. Observe that there is a forgetful functor \func{U}{{\rm Sh}_{\ct{E}}(\ct{C})}{\ct{E}/C_0}, and call a morphism \func{f}{B}{A} of sheaves \emph{pointwise small}, when $Uf$ is. To show that these morphisms form a class of small maps, we make two additional assumptions. First of all, we assume the Exponentiation Axiom in the ``metatheory'' 
$(\ct{E}, \smallmap{S})$:
\begin{description}
\item[($\Pi$S)] For any small map \func{f}{B}{A}, the functor $\Pi_f: \ct{E}/B \to \ct{E}/A$ preserves small objects.
\end{description}
Furthermore, we also assume that our site has a basis, meaning the following: for any $a \in C_0$ there is a \emph{small} collection of covering sieves ${\rm BCov}(a)$ such that
\[ S \in {\rm Cov}(a) \Leftrightarrow \exists R \in {\rm BCov}(a): R \subseteq S. \]
Note that the relation $S \in {\rm Cov}(a)$ is bounded, when the site has a basis.

\begin{theo}{stabundersheaves}
Let $(\ct{E}, \smallmap{S})$ be a category with small maps, and let \ct{C} be an internal site with a basis. If the class \smallmap{S} satisfies {\bf ($\Pi$S)}, then ${\rm Sh}_{\ct{E}}(\ct{C})$ with the class of pointwise small maps is again a category with small maps satisfying {\bf ($\Pi$S)}. Furthermore, all the axioms that we have introduced, {\bf (NS)}, {\bf (PS)}, {\bf (F)} and {\bf (M)}, are stable in the sense that each of these holds in sheaves, whenever it holds in the original category.
\end{theo}

\section{Predicative realisability}

In this Section we outline how the construction of \cite{joyalmoerdijk95} of a class of small maps in Hyland's effective topos (as in Example 5.3), can be mimicked in the context of a category with small maps $(\ct{E}, \smallmap{S})$ as introduced in Section 3. Our construction is inspired by the fact that the effective topos arises as the exact completion of the category of assemblies, as in \cite{carbonifreydscedrov88}. 

Let us start with a category with small maps $(\ct{E}, \smallmap{S})$ satisfying {\bf (NS)} (so the nno in \ct{E} is small). The first observation is that we can internalise enough recursion theory in \ct{E} for doing realisability. In fact, enough can already be formalised in Heyting Arithmetic {\bf HA}, so certainly in a category with small maps. We then define the category of assemblies, as follows. An \emph{assembly} consists of an object $A$ in \ct{E} together with a surjective relation $\alpha \subseteq \NN \times A$. For pairs $(n, a)$ belonging to this relation, we write $n \in \alpha(a)$, which we pronounce as ``$n$ realises (the existence of) $a$''; surjectivity of the relation then means that every $a \in A$ has at least one realiser. A morphism $f$ of assemblies from $(B, \beta)$ to $(A, \alpha)$ is given by a morphism \func{f}{B}{A} in \ct{E} for which the internal logic of \ct{E} verifies that:
\begin{quote}
there is a natural number $r$ such that for all $b \in B$ and $n \in \beta(b)$, the Kleene application $r \cdot n$ is defined, and realises $f(b)$ (i.e.~$r \cdot n \in \alpha(fb)$).
\end{quote}
One can now prove that the category of assemblies $\ct{E}[\Asm]$ relative to \ct{E} is a Heyting category with stable and disjoint sums (see \cite{hyland82}, where the assemblies occur as the $\lnot\lnot$-separated objects in the effective topos).

In order to describe the relevant exact completion of this category of assemblies, we first outline a construction. Consider two assemblies $(B, \beta)$ and $(A, \alpha)$ and a morphism \func{f}{B}{A}, not necessarily a morphism of assemblies. Then this defines a morphism of assemblies $(B, \beta[f]) \rTo (A, \alpha)$ by declaring that $n \in \beta[f](b)$, whenever $n$ codes a pair $< n_0, n_1 >$ such that $n_0 \in \alpha(fb)$ and $n_1 \in \beta(b)$. In case $f$ belongs to \smallmap{S} and $\beta$ is a bounded relation, a morphism of this form will be called a \emph{standard display map} relative to \smallmap{S} (this notion was pointed out to us by Thomas Streicher). A \emph{display map} is a morphism that can be written as an isomorphism followed by a standard display map. These display maps do \emph{not} satisfy the axioms for a class of small maps; in particular, they are not closed under Descent and Quotients. Another problem is that the category of assemblies is not exact, not even in the more limited sense of being bounded exact.

Both problems can be solved by appealing to \reftheo{consofBE}. Or, to be more precise, they can be solved by constructing an exact completion for categories with a class of display maps, resulting in categories with small maps satisfying {\bf (BE)} (how this is to be done will be shown in \cite{bergmoerdijk07}). Recall that the small maps in the exact completion are precisely those $g$ that fit into a quasi-pullback diagram
\diag{{\bf y}D \ar[d]_{{\bf y}f} \ar@{->>}[r] & B \ar[d]^g \\
{\bf y}C \ar@{->>}[r] & A, }
where $f$ is a small map in the original category. Therefore it is to be expected that the class of small maps in the exact completion of a category with display maps satisfies Descent and Quotients even when the class of display maps in the original category from which it is defined, does not satisfy these axioms. In fact, as it turns out, the display maps between assemblies have enough structure for the maps $g$ in the exact completion of assemblies that fit into a square as above with $f$ a display map, to form a class of small maps. In this way, both problems with the category of assemblies can be solved at the same time by moving to the exact completion. Therefore we define the realisability category $(\ct{E}[\Eff], \smallmap{S}[\Eff]])$ to be this exact completion of the pair $(\ct{E}[\Asm],\smallmap{D})$, where \smallmap{D} is the class of display maps in the category of assemblies.

\begin{theo}{stabunderrealisability}
If $(\ct{E}, \smallmap{S})$ is a category with small maps satisfying {\bf (NS)}, then so is $(\ct{E}[\Eff], \smallmap{S}[\Eff]])$. Furthermore, all the axioms that we have introduced, {\bf ($\Pi$S)}, {\bf (PS)}, {\bf (F)} and {\bf (M)}, are stable in the sense that each of these holds in the realisablity category, whenever it holds in the original category.
\end{theo}

The initial ZF-algebra in the realisability category should be considered as a suitable internal version of McCarty's realisability model \cite{mccarty84} (see also \cite{kouwenhovenvanoosten05}), which in our abstract approach is also defined for predicative theories like {\bf CZF} (compare \cite{rathjen06}).

\bibliographystyle{plain} \bibliography{ast}

\begin{thebibliography}{10}

\bibitem{aczel78}
P.~Aczel.
\newblock The type theoretic interpretation of constructive set theory.
\newblock In {\em Logic Colloquium '77 (Proc. Conf., Wroc\l aw, 1977)},
  volume~96 of {\em Stud. Logic Foundations Math.}, pages 55--66. North-Holland
  Publishing Co., Amsterdam, 1978.

\bibitem{aczelrathjen01}
P.~Aczel and M.~Rathjen.
\newblock Notes on constructive set theory.
\newblock Technical Report No. 40, Institut Mittag-Leffler, 2000/2001.

\bibitem{awodeyetal04}
S.~Awodey, C.~Butz, A.K. Simpson, and T.~Streicher.
\newblock Relating topos theory and set theory via categories of classes.
\newblock Available from http://www.phil.cmu.edu/projects/ast/, June 2003.

\bibitem{awodeyforssell05}
S.~Awodey and H.~Forssell.
\newblock Algebraic models of intuitionistic theories of sets and classes.
\newblock {\em Theory Appl. Categ.}, 15:No. 5, 147--163 (electronic), 2005/06.

\bibitem{awodeywarren05}
S.~Awodey and M.A. Warren.
\newblock Predicative algebraic set theory.
\newblock {\em Theory Appl. Categ.}, 15:1, 1--39 (electronic), 2005/06.

\bibitem{berg05}
B.~van~den Berg.
\newblock Inductive types and exact completion.
\newblock {\em Ann. Pure Appl. Logic}, 134:95--121, 2005.

\bibitem{berg06}
B.~van~den Berg.
\newblock {\em Predicative topos theory and models for constructive set
  theory}.
\newblock PhD thesis, University of Utrecht, 2006.

\bibitem{berg05a}
B.~van~den Berg.
\newblock Sheaves for predicative toposes.
\newblock Accepted for publication in \emph{Arch. Math. Logic}. Available from
  arXiv: math.LO/0507480, July 2005.

\bibitem{bergdemarchi06}
B.~van~den Berg and F.~De~Marchi.
\newblock Models of non-well-founded sets via an indexed final coalgebra
  theorem.
\newblock Accepted for publication in \emph{J. Symbolic Logic}. Available from
  arXiv: math.LO/0508531, Jan 2006.

\bibitem{bergmoerdijk07}
B.~van~den Berg and I.~Moerdijk.
\newblock Predicative algebraic set theory, {I}: exact completion, {II}: sheaf
  models, {III}: realisability models.
\newblock Under construction, 2007.

\bibitem{blassscedrov89}
A.R. Blass and A.~Scedrov.
\newblock Freyd's models for the independence of the axiom of choice.
\newblock {\em Mem. Am. Math. Soc.}, 79(404), 1989.

\bibitem{blassscedrov92}
A.R. Blass and A.~Scedrov.
\newblock Complete topoi representing models of set theory.
\newblock {\em Ann. Pure Appl. Logic}, 57(1):1--26, 1992.

\bibitem{butz03}
C.~Butz.
\newblock Bernays-{G}\"odel type theory.
\newblock {\em J. Pure Appl. Algebra}, 178(1):1--23, 2003.

\bibitem{carboni95}
A.~Carboni.
\newblock Some free constructions in realizability and proof theory.
\newblock {\em J. Pure Appl. Algebra}, 103:117--148, 1995.

\bibitem{carboniceliamagno82}
A.~Carboni and R.~Celia~Magno.
\newblock The free exact category on a left exact one.
\newblock {\em J. Austral. Math. Soc.}, 33:295--301, 1982.

\bibitem{carbonifreydscedrov88}
A.~Carboni, P.J. Freyd, and A.~Scedrov.
\newblock A categorical approach to realizability and polymorphic types.
\newblock In {\em Mathematical foundations of programming language semantics
  (New Orleans, LA, 1987)}, volume 298 of {\em Lecture Notes in Comput. Sci.},
  pages 23--42. Springer Verlag, Berlin, 1988.

\bibitem{fourman80}
M.P. Fourman.
\newblock Sheaf models for set theory.
\newblock {\em J. Pure Appl. Algebra}, 19:91--101, 1980.

\bibitem{freyd80}
P.J. Freyd.
\newblock The axiom of choice.
\newblock {\em J. Pure Appl. Algebra}, 19:103--125, 1980.

\bibitem{gambino05}
N.~Gambino.
\newblock Presheaf models for constructive set theories.
\newblock In {\em From sets and types to topology and analysis}, volume~48 of
  {\em Oxford Logic Guides}, pages 62--77. Oxford University Press, Oxford,
  2005.

\bibitem{gambino06}
N.~Gambino.
\newblock Heyting-valued interpretations for constructive set theory.
\newblock {\em Ann. Pure Appl. Logic}, 137(1-3):164--188, 2006.

\bibitem{gambinohyland04}
N.~Gambino and J.M.E. Hyland.
\newblock Wellfounded trees and dependent polynomial functors.
\newblock In {\em Types for proofs and programs}, volume 3085 of {\em Lecture
  Notes in Comput. Sci.}, pages 210--225. Springer-Verlag, Berlin, 2004.

\bibitem{grayson83}
R.J. Grayson.
\newblock Forcing in intuitionistic systems without power-set.
\newblock {\em J. Symbolic Logic}, 48(3):670--682, 1983.

\bibitem{hyland82}
J.M.E. Hyland.
\newblock The effective topos.
\newblock In {\em The L.E.J. Brouwer Centenary Symposium (Noordwijkerhout,
  1981)}, volume 110 of {\em Stud. Logic Foundations Math.}, pages 165--216.
  North-Holland Publishing Co., Amsterdam, 1982.

\bibitem{hylandrobinsonrosolini90}
J.M.E. Hyland, E.P. Robinson, and G.~Rosolini.
\newblock The discrete objects in the effective topos.
\newblock {\em Proc. London Math. Soc. (3)}, 60(1):1--36, 1990.

\bibitem{johnstone02a}
P.T. Johnstone.
\newblock {\em Sketches of an elephant: a topos theory compendium. {V}olume 1},
  volume~43 of {\em Oxf. Logic Guides}.
\newblock Oxford University Press, New York, 2002.

\bibitem{johnstone02b}
P.T. Johnstone.
\newblock {\em Sketches of an elephant: a topos theory compendium. {V}olume 2},
  volume~44 of {\em Oxf. Logic Guides}.
\newblock Oxford University Press, Oxford, 2002.

\bibitem{joyalmoerdijk95}
A.~Joyal and I.~Moerdijk.
\newblock {\em Algebraic set theory}, volume 220 of {\em London Mathematical
  Society Lecture Note Series}.
\newblock Cambridge University Press, Cambridge, 1995.

\bibitem{kouwenhovenvanoosten05}
C.~Kouwenhoven-Gentil and J.~van Oosten.
\newblock Algebraic set theory and the effective topos.
\newblock {\em J. Symbolic Logic}, 70(3):879--890, 2005.

\bibitem{lambekscott86}
J.~Lambek and P.J. Scott.
\newblock {\em Introduction to higher order categorical logic}, volume~7 of
  {\em Camb. Stud. Adv. Math.}
\newblock Cambridge University Press, Cambridge, 1986.

\bibitem{lubarsky06}
R.S. Lubarsky.
\newblock {CZF} and {S}econd {O}rder {A}rithmetic.
\newblock {\em Ann. Pure Appl. Logic}, 141(1-2):29--34, 2006.

\bibitem{mccarty84}
D.C. McCarty.
\newblock {\em Realizability and recursive mathematics}.
\newblock PhD thesis, Oxford University, 1984.

\bibitem{moerdijkpalmgren00}
I.~Moerdijk and E.~Palmgren.
\newblock Wellfounded trees in categories.
\newblock {\em Ann. Pure Appl. Logic}, 104(1-3):189--218, 2000.

\bibitem{moerdijkpalmgren02}
I.~Moerdijk and E.~Palmgren.
\newblock Type theories, toposes and constructive set theory: predicative
  aspects of {AST}.
\newblock {\em Ann. Pure Appl. Logic}, 114(1-3):155--201, 2002.

\bibitem{rathjen06}
M.~Rathjen.
\newblock Realizability for constructive {Z}ermelo-{F}raenkel set theory.
\newblock In {\em Logic Colloquium '03}, volume~24 of {\em Lect. Notes Log.},
  pages 282--314. Assoc. Symbol. Logic, La Jolla, CA, 2006.

\bibitem{simpson99}
A.K. Simpson.
\newblock Elementary axioms for categories of classes (extended abstract).
\newblock In {\em 14th Symposium on Logic in Computer Science (Trento, 1999)},
  pages 77--85. IEEE Computer Soc., Los Alamitos, CA, 1999.

\bibitem{simpson05}
A.K. Simpson.
\newblock Constructive set theories and their category-theoretic models.
\newblock In {\em From sets and types to topology and analysis}, volume~48 of
  {\em Oxford Logic Guides}, pages 41--61. Oxford University Press, Oxford,
  2005.

\bibitem{streicher05}
T.~Streicher.
\newblock Realizability models for {CZF}+ $\lnot$ {P}ow.
\newblock Unpublished note, March 2005.

\end{thebibliography}

\end{document}